\documentclass[a4paper,11pt]{amsart}
\usepackage{amssymb}

\usepackage{amsmath}

\newcommand{\F}{\mathbb{F}}

\DeclareMathOperator{\ord}{ord}

\newtheorem{dummy}{Dummy}

\numberwithin{dummy}{section}
\numberwithin{equation}{section}

\newtheorem{lemma}[dummy]{Lemma}
\newtheorem{theorem}[dummy]{Theorem}

\newtheorem{cor}[dummy]{Corollary}

\theoremstyle{definition}

\theoremstyle{remark}
\newtheorem{rem}[dummy]{Remark}

\begin{document}
\bibliographystyle{amsalpha}

\date{\today}
\author{S.~Mattarei}
\email{mattarei@science.unitn.it}
\urladdr{http://www-math.science.unitn.it/\~{ }mattarei/}
\address{Dipartimento di Matematica\\
  Universit\`a degli Studi di Trento\\
  via Sommarive 14\\
  I-38050 Povo (Trento)\\
  Italy}
\title[The orders of nonsingular derivations]%
{A sufficient condition for a number to be the order of a nonsingular derivation of a Lie algebra}

\begin{abstract}
A study of the set $\mathcal{N}_p$ of positive integers which occur as
orders of nonsingular derivations of finite-dimensional non-nilpotent Lie algebras
of characteristic $p>0$ was initiated by Shalev and continued by the present author.
The main goal of this paper is to produce more elements of $\mathcal{N}_p$.
Our main result shows that any divisor $n$ of $q-1$, where $q$ is a power of $p$, such that
$n\ge (p-1)^{1/p} (q-1)^{1-1/(2p)}$, necessarily belongs to $\mathcal{N}_p$.
This extends its special case for $p=2$ which was proved in a previous paper by a different method.
\end{abstract}



\maketitle

\thispagestyle{empty}

\section{Introduction}\label{Intro}

Finite-dimensional Lie algebras which admit a nonsingular (that is, injective) derivation
play a role in various investigations.
Some of those are discussed in the Introduction of~\cite{Mat:nonsing-der-2}, of which this paper is a continuation.
We briefly recall here only the essential facts relevant to our present study and refer to~\cite{Mat:nonsing-der-2}
and its predecessor~\cite{Mat:nonsing-der} for more details.

According to a classical result of Jacobson~\cite[p.~54]{Jac:Lie_algebras},
in characteristic zero only nilpotent Lie algebras can have nonsingular derivations.
In positive characteristic, where even certain simple Lie algebras can admit nonsingular derivations,
the same argument would be inconclusive, but still imposes a strong restriction
of the eigenvalues (assumed in the ground field) of a nonsingular derivation of a non-nilpotent Lie algebra.
In particular, if the derivation has finite order $n$, as is relevant to various studies,
this restriction entails an interesting necessary condition on $n$, noted by Shalev in~\cite{Sha:nonsing-der}.
The condition was shown to be sufficient as well in~\cite{Mat:nonsing-der}.
We recall both implications as Theorem~\ref{thm:N_p} in the next section.

More generally, in his paper~\cite{Sha:nonsing-der} Shalev suggested and initiated a study of
the set $\mathcal{N}_p$ of positive integers which occur as the orders of nonsingular
derivations of finite-dimensional non-nilpotent Lie algebras of prime characteristic $p$.
Theorem~\ref{thm:N_p} translates this problem into one entirely formulated in terms of finite fields.
Therefore, no Lie algebra argument will be used in this paper.
It is easy to see that $\mathcal{N}_p$ is closed with respect to taking multiples, and that
a positive integer $n$ belongs to $\mathcal{N}_p$ if and only if its $p'$-part does.
Thus, one may restrict one's attention to numbers in $\mathcal{N}_p$ which are prime to $p$.
Even after this restriction, rather trivial elements of $\mathcal{N}_p$ are those of the form
$p^k-1$ for some $k\ge 2$, as will be clear from Theorem~\ref{thm:N_p}.
We will conveniently call {\em nontrivial} elements of $\mathcal{N}_p$
those numbers in $\mathcal{N}_p$ which are prime to $p$
and are not multiples of any $p^k-1$ with $k\ge 2$.
Shalev proved in~\cite{Sha:nonsing-der} that no nontrivial element of $\mathcal{N}_p$
is smaller than $p^2$.
This was extended in~\cite[Lemma~3.2]{Mat:nonsing-der},
to conclude that no nontrivial element of $\mathcal{N}_p$ is smaller than $p^3$,
except for $(3^3-1)/2=13$ when $p=3$.
(This exception has an analogue for all odd primes, see the next paragraph.)
In fact, we do not know of any nontrivial element of $\mathcal{N}_p$
which is smaller than $p^4$, except for $13$ when $p=3$.

In the opposite direction, one nontrivial element of $\mathcal{N}_p$ is $(p^p-1)/(p-1)$, for odd $p$,
as noted by Shalev in~\cite[Example~2.6]{Sha:nonsing-der}.
For $p=2$ many numbers in $\mathcal{N}_2$ were disclosed in~\cite{Mat:nonsing-der-2}.
Apart from the special series of numbers of the form $n=(2^{3s}-1)/(2^s-1)$, of which the case $s=3$ (whence $n=73$) was
already noted by Shalev in~\cite{Sha:nonsing-der},
we proved a result guaranteeing that all divisors of $q-1$, where $q$ is a power of $2$,
which are large enough in an appropriate sense belong to $\mathcal{N}_2$.
Explicitly, a sufficient condition for a divisor $n$ of $q-1$ to belong to $\mathcal{N}_2$ was found to be that $n\ge (q-1)^{3/4}$.
However, the arguments used in~\cite{Mat:nonsing-der-2},
based on the character theory of a certain group, were limited to the case of the prime $2$.

In this paper we extend that result to an arbitrary prime $p$.
We prove in Corollary~\ref{cor:n_large} that a divisor $n$ of $q-1$, where $q$ is a power of $p$,
belongs to $\mathcal{N}_p$ provided it satisfies the inequality
$n\ge (p-1)^{1/p} (q-1)^{1-1/(2p)}$.
This is a weaker and simplified form of a more precise sufficient condition
for a certain system of equations having solutions over the finite field $\F_q$.
We prove that in Section~\ref{sec:system} by means of standard character sum estimates.
We sketch a less elementary but shorter proof in Remark~\ref{rem:sophisticated}.

In Section~\ref{sec:comments} we have collected several remarks on the set $\mathcal{N}_p$.
In particular, we discuss some consequences of our main result,
present the outcome of some computer calculations,
discuss the density of the set of integers $\mathcal{N}_p$
(following a suggestion of a referee)
and a notion of relative size of its elements.

\section{Large divisors of $q-1$ belong to $\mathcal{N}_p$}\label{sec:large}

As mentioned in the Introduction,
$\mathcal{N}_p$ denotes the set of positive integers $n$
for which there exists finite-dimensional non-nilpotent Lie algebra $L$, over a field of characteristic $p$,
which admits a nonsingular derivation of order $n$.
We recall from~\cite[Corollary~2.3]{Mat:nonsing-der} the essential part
of a characterization of the elements of $\mathcal{N}_p$ which are prime to $p$.

\begin{theorem}\label{thm:N_p}
A positive integer $n$ prime to $p$ belongs to $\mathcal{N}_p$ if and only if
there exists an element $\xi$ of the algebraic closure $\bar\F_p$ of $\F_p$,
such that
$(\xi+\lambda)^n=1$ for all $\lambda\in\F_p$.
\end{theorem}

This condition is trivially satisfied by numbers $n$ of the form $p^k-1$ with $k\ge 2$,
and hence those numbers belong to $\mathcal{N}_p$, together with their multiples.
As anticipated in the Introduction, we call {\em trivial} those elements of $\mathcal{N}_p$,
and {\em nontrivial} the others.

In general, for any $n$ prime to $p$ there is a power $q$ of $p$ such that $n$ divides $q-1$.
For example, we may take $q=p^k$, where $k$ is the multiplicative order of $p$ modulo $n$.
Then the condition $n\in\mathcal{N}_p$ is equivalent to the fact that there exists an element $\xi$
of the finite field $\F_q$, such that $\xi,\xi+1,\ldots,\xi+p-1$ are nonzero $d$th powers in $\F_q$,
where $d=(q-1)/n$.
The following result provides an estimate for the number of such elements $\xi$,
in a more general setting.

\begin{theorem}\label{thm:sequence}
Let $d$ be a divisor of $q-1$
and let $0<r\le p$.
Let $M$ be the number of elements $\xi$ of $\F_q$ such that
$\xi,\xi+1,\ldots,\xi+r-1$ are nonzero $d$th powers in $\F_q$.
Let $M_0$ be the number of elements $\xi$ of $\F_p$ such that
$\xi,\xi+1,\ldots,\xi+r-1$ include $0$ and are $d$th powers in $\F_q$.
Then
\[
\left|M+\frac{M_0+1}{d}-\frac{q+1}{d^r}\right|
\le
\left(r-1-\frac{r+1}{d}+\frac{2}{d^r}\right)
\sqrt{q}.
\]
\end{theorem}

Since $0\le M_0\le r$ we deduce the bound
\[
\left|M+\frac{r+1}{2d}-\frac{q+1}{d^r}\right|
\le
\left(r-1-\frac{r+1}{d}+\frac{2}{d^r}\right)
\sqrt{q}+\frac{r+1}{2d}.
\]
which does not involve $M_0$.
Bounds of this type result from standard calculations with character sums, see~\cite[Exercises~5.65 and~5.66]{LN}
or~\cite[pp.~246-247]{Stepanov:arithmetic}.
Their simplest application is that, given $d$ and $r$, for all primes $p$ large enough there exists a
sequence of $r$ consecutive integers which are $d$th power residues modulo $p$.
However, we are unable to quote from the literature a bound which is as sharp as that given in Theorem~\ref{thm:sequence}
(see the discussion in Remark~\ref{rem:weak}), and hence we provide a proof in Section~\ref{sec:system}.

Here we need the special case of Theorem~\ref{thm:sequence} where $r=p$.
Then the lower bound for $M$ reads
\begin{equation}\label{eq:lower}
M\ge\frac{q+1}{d^p}-\frac{M_0+1}{d}
-
\left(p-1-\frac{p+1}{d}+\frac{2}{d^{p}}\right)
\sqrt{q},
\end{equation}
and $M_0$ can only be $p$ or $0$, according as $n$ is a multiple of $p-1$ or not.
Then we know that $n$ belongs to $\mathcal{N}_p$ exactly when $M>0$, where $d=(q-1)/n$.
Thus, a sufficient condition for $n\in\mathcal{N}_p$ is that the right-hand side of inequality~\eqref{eq:lower}
be strictly positive.
After a simple calculation this yields the following result.

\begin{theorem}\label{thm:n_large}
Let $q$ be a power of $p$ and let $n$ be a divisor of $q-1$.
Then $n\in\mathcal{N}_p$ provided $d=(q-1)/n$ satisfies
\[
q-((pd-p-d-1)d^{p-1}+2)\sqrt{q}-(p+1)d^{p-1}+1>0.
\]
If $n$ is not a multiple of $(p-1)$, then the slightly weaker condition
\[
q-((pd-p-d-1)d^{p-1}+2)\sqrt{q}-d^{p-1}+1>0
\]
suffices.
\end{theorem}

\begin{rem}\label{rem:Weil}
For $p=2$, where $n$ is, necessarily, a multiple of $p-1$, the sufficient condition
for $n\in\mathcal{N}_p$ given in Theorem~\ref{thm:n_large} reads
$q-(d-1)(d-2)\sqrt{q}-3d+1>0$.
Once expressed in terms of $n=(q-1)/d$, the condition becomes
$n^2+3(\sqrt{q}+1)n-\sqrt{q}(\sqrt{q}+1)^2>0$
or, equivalently,
$n>(\sqrt{4\sqrt{q}+9}-3)(\sqrt{q}+1)/2$.
This is slightly weaker than the sufficient condition $n^4>(q-n)^3$ given
in~\cite[Theorem~3.1]{Mat:nonsing-der-2}.
The reason is the following.
As will be clear after section~\ref{sec:system}, especially Remark~\ref{rem:weak}, when $p=2$
the sufficient condition of Theorem~\ref{thm:n_large} ultimately depends on Weil's bound
$|N-d-q-1|\le(d-1)(d-2)\sqrt{q}$ for the number $N$ of affine points of the Fermat curve $y_2^d-y_1^d=1$ over $\F_q$.
One can see that the proof of~\cite[Theorem~3.1]{Mat:nonsing-der-2} establishes and then uses
a weaker bound than Weil's, with an error term close to
$(d^2-\frac{3}{2}d)\sqrt{q}$
rather than $(d^2-3d+2)\sqrt{q}$.
\end{rem}

A slightly weaker but more manageable form of the sufficient conditions given in Theorem~\ref{thm:n_large}
is the following.

\begin{cor}\label{cor:n_large}
Let $q$ be a power of $p$ and let $n$ be a divisor of $q-1$ such that
\[
n\ge (p-1)^{1/p} (q-1)^{1-1/(2p)}.
\]
Then $n\in\mathcal{N}_p$.
\end{cor}

Note that the factor $(p-1)^{1/p}$ is always less than $1.32$ and rapidly tends to $1$ as $p$ tends to infinity.
When $p=2$ the condition in Corollary~\ref{cor:n_large} reads $n\ge(q-1)^{3/4}$,
which is only slightly stronger than the condition $n^4>(q-n)^3$ of~\cite[Theorem~3.1]{Mat:nonsing-der-2}.

\begin{proof}
The former (and stronger) inequality in Theorem~\ref{thm:n_large} can be equivalently written as
\[
\bigl(\sqrt{q}-(p-1)d^p\bigr)
\bigl(\sqrt{q}+(p+1)d^{p-1}-2\bigr)
+
\bigl((p-1)d^p-1\bigr)
\bigl((p+1)d^{p-1}-2\bigr)
>1.
\]
Temporarily viewing $\sqrt{q}$ as a real variable,
the inequality is satisfied when $\sqrt{q}=(p-1)d^p$,
except when $p-1=d=1$.
But in that case the conclusion of Corollary~\ref{cor:n_large} holds trivially.
Consequently, the inequality holds whenever $\sqrt{q}\ge (p-1)d^p$.
In particular, it holds whenever $q-1\ge (p-1)^2d^{2p}$,
which is equivalent with the stated hypothesis when written in terms of $n=(q-1)/d$.
\end{proof}

\begin{rem}\label{rem:asymptotic}
Using the form of the inequality used in the proof of Corollary~\ref{cor:n_large},
one can easily see that the sufficient condition for $n\in\mathcal{N}_p$
given in Theorem~\ref{thm:n_large} is asymptotic to the simpler one given in Corollary~\ref{cor:n_large},
in the sense that
\[
\lim_{q\to\infty}(p-1)^{1/p} (q-1)^{1-1/(2p)}/f(q)=1,
\]
where $n>f(q)$ is an explicit form of the condition given in the former.
\end{rem}

\section{Comments, calculations, further questions}\label{sec:comments}

\subsection{Existence of proper divisors of $p^k-1$ in $\mathcal{N}_p$}
We discuss in which respect our main result produces nontrivial elements of $\mathcal{N}_p$.
The following is an essentially equivalent formulation of Corollary~\ref{cor:n_large}
in terms of $d=(q-1)/n$ in place of $n$.

\begin{cor}\label{cor:k_large}
Let $p$ be a prime and let $d$ be a positive integer prime to $p$.
If $k$ is a positive multiple of the order of $p$ modulo $d$
then $(p^k-1)/d\in\mathcal{N}_p$
provided
$k\ge 2+2p\,\log d/\log p$.
\end{cor}

\begin{proof}
Setting $q=p^k$ and in terms of $d=(q-1)/n$, the sufficient condition of Corollary~\ref{cor:n_large} becomes
$p^k\ge(p-1)^2d^{2p}+1$, but the proof of Corollary~\ref{cor:n_large} shows that the summand $1$
can be discarded.
Our present hypothesis
$k\ge 2+2p\,\log d/\log p$
is only slightly stronger than that.
\end{proof}

By taking $d=2$ in Corollary~\ref{cor:k_large} we see that, for every odd prime $p$ and every integer $k\ge 2+p\,\log 4/\log p$,
there is at least one proper divisor of $p^k-1$ which belongs to $\mathcal{N}_p$,
namely, $(p^k-1)/2$.
This statement is actually nontrivial only when $k$ is prime,
because otherwise $p^k-1$ has proper divisors of the form $p^s-1$ with $s>1$,
which are trivial elements of $\mathcal{N}_p$.
Incidentally, note that
the simplified condition of Corollary~\ref{cor:k_large}
(as well as that of Corollary~\ref{cor:n_large})
is notably weaker than the more precise Theorem~\ref{thm:n_large}
for small $p$ and $k$.
For example, Corollary~\ref{cor:k_large} implies that $(3^k-1)/2\in\mathcal{N}_3$ for $k\ge 6$,
while the inequalities in Theorem~\ref{thm:n_large} show that this is the case for $k=3,4,5$ as well.

\subsection{Varying the characteristic}
It is also interesting to look at Corollary~\ref{cor:k_large},
or to the more precise Theorem~\ref{thm:n_large} when needed, 
from a different perspective, thinking of $k$ as assigned and varying the prime $p$.
The smallest value of $k$ which is of interest here is $k=5$.
In fact, according to~\cite[Corollary~3.4]{Mat:nonsing-der},
no proper divisor of $p^3-1$ belongs to $\mathcal{N}_p$, with the only exception that $(3^3-1)/2=13\in\mathcal{N}_3$.
Moreover, when $k=4$ and $p>2$
the number $(p^4-1)/2$ is a multiple of $p^2-1$ and, hence, is a trivial element of $\mathcal{N}_p$.
When $k=5$, Theorem~\ref{thm:n_large} implies that $(p^5-1)/2\in\mathcal{N}_p$ for $p=3,5$, as mentioned above.
(As reported in~\cite[Example~4.1]{Mat:nonsing-der}, direct calculations show that $(p^5-1)/2\in\mathcal{N}_p$
for $p=7,11$ as well, but not for $p=13$.)
Similarly, Theorem~\ref{thm:n_large} implies that $(p^7-1)/2\in\mathcal{N}_p$ for $p=3,5,7$.

In this respect we should note that, for a fixed prime $k$,
there can only be finitely many primes $p$
such that $p^k-1$ has a proper divisor in $\mathcal{N}_p$.
This follows from a result of H.~Davenport~\cite[Theorem~1]{Dav}:
given $k>1$ (not necessarily prime), if the prime $p$ is sufficiently large
(depending only on $k$)
and $\F_{p^k}=\F_p(\xi)$, then there exists $\lambda\in\F_p$
such that $\xi+\lambda$ is a primitive element of $\F_{p^k}$.
Under our present assumption that $k$ is prime, any element $\xi\in\F_{p^k}$ chosen as in Theorem~\ref{thm:N_p}
satisfies $\F_{p^k}=\F_p(\xi)$, and then Davenport's theorem implies that $p^k-1$ divides $n$
if $p$ is sufficiently large.

\subsection{Computer calculations}\label{subsec:calculations}
Since $\mathcal{N}_p$ is closed with respect to taking multiples, the following definition is convenient:
call {\em minimal} any number in $\mathcal{N}_p$ which has no proper divisor in $\mathcal{N}_p$.
(In the terminology of~\cite[Chapter V]{HalRot},
the minimal elements of $\mathcal{N}_p$ form the {\em primitive generating sequence} of $\mathcal{N}_p$.)
A computer search has shown that the minimal elements of $\mathcal{N}_2$ below $200000$ are
\begin{align*}
&
3=2^2-1,\quad
7=2^3-1,\quad
31=2^5-1,\quad
73=(2^9-1)/7,\quad
85=(2^8-1)/3,
\\&
127=2^7-1,\quad
2047=2^{11}-1,\quad
3133=(2^{24}-1)/5355,\quad
4369=(2^{16}-1)/15,
\\&
8191=2^{13}-1,\quad
11275=(2^{20}-1)/93,\quad
49981=(2^{30}-1)/21483,
\\&
60787=(2^{22}-1)/69,\quad
76627=(2^{36}-1)/896805,\quad
121369=(2^{39}-1)/4529623,
\\&
131071=2^{17}-1,\quad
140911=(2^{28}-1)/1905,\quad
178481=(2^{23}-1)/47.
\end{align*}
Of the nontrivial elements in this list,
only $85$ and $4369$ are explained by Corollary~\ref{cor:n_large}.
For other elements in the list, Corollary~\ref{cor:n_large} (or Theorem~\ref{thm:n_large}) is only strong enough to show
that certain of their multiples, still within the range considered, are nontrivial elements of $\mathcal{N}_p$.
As an example, this is the case for $11275\cdot 3=(2^{20}-1)/31$.
Another fact which follows from inspection of the table,
together with the observation that $2^{19}-1$ is a prime,
is that $2^{23}-1$ is the smallest element of $\mathcal{N}_2$, of the form $2^k-1$ with $k$ prime,
which is not minimal.

We have carried out similar calculations for $p=3$.
They have shown that the minimal elements of $\mathcal{N}_3$ below $100000$ are
\begin{align*}
&
8=3^2-1,\quad
13=(3^3-1)/2,\quad
121=(3^5-1)/2,
\\
&
1093=(3^7-1)/2,\quad
88573=(3^{11}-1)/2,
\end{align*}
all of which are predicted by Theorem~\ref{thm:n_large}, as discussed above.
The smallest prime $k$ for which Theorem~\ref{thm:n_large} produces a proper divisor of $3^k-1$ in $\mathcal{N}_3$
different from $(3^k-1)/2$
is $23$, namely, we have $(3^{23}-1)/47\in\mathcal{N}_3$.
Note that such a number would be far too large for a direct verification that it belongs to $\mathcal{N}_3$
based on the characterization given in Theorem~\ref{thm:N_p}.
Because of the computational complexity of an exhaustive search we were not able
to produce any element of $\mathcal{N}_3$ which is not predicted by Theorem~\ref{thm:n_large}.

For each prime $p$ larger than $3$ we know of essentially only one element of $\mathcal{N}_p$
which is not within the range where Theorem~\ref{thm:n_large} applies,
namely, the number $(p^p-1)/(p-1)$ noted by Shalev in~\cite[Example~2.6]{Sha:nonsing-der}.
This is also the smallest element of $\mathcal{N}_p$ which we know of for a generic prime $p>3$.
(The element $(p^k-1)/2$ produced by Corollary~\ref{cor:k_large} is larger than that.)

\subsection{Density of $\mathcal{N}_p$}
A referee has suggested to look at the density of $\mathcal{N}_p$.
It is not clear whether $\mathcal{N}_p$ possesses a {\em natural density}
$
\lim_{m\to\infty}|\{n\in\mathcal{N}_p\mid n\le m\}|/m.
$
However, since $\mathcal{N}_p$ is closed with respect to taking multiples,
a result of Davenport and Erd\"os~\cite[Chapter~V, Theorem~12]{HalRot} guarantees that $\mathcal{N}_p$ possesses a
{\em logarithmic density}
\[
\delta(\mathcal{N}_p)=
\lim_{m\to\infty}\frac{1}{\log m}\sum_{n\in\mathcal{N}_p,\ n\le m}\frac{1}{n},
\]
and that $\delta(\mathcal{N}_p)$ coincides with the {\em lower asymptotic density,}
$
\liminf_{m\to\infty}|\{n\in\mathcal{N}_p\mid n\le m\}|/m.
$

It would be interesting to know how much $\delta(\mathcal{N}_p)$
exceeds $\delta(\mathcal{T}_p)$, where $\mathcal{T}_p$
consists of the trivial elements of $\mathcal{N}_p$, that is, of all multiples of numbers of the form $p^k-1$ for some $k\ge 2$.
Since $\mathcal{T}_p$ coincides with the set of multiples of numbers of the form
$p^k-1$ with $k$ prime, and because $p^k-1$ and $p^{k'}-1$ have greatest common divisor $p-1$
for different primes $k$ and $k'$, one easily sees that $\mathcal{T}_p$
possesses a natural density, whose value equals
\[
\delta(\mathcal{T}_p)=
\frac{1}{p-1}\left(1-\prod_{k\text{ prime}}
\left(1-\frac{p-1}{p^k-1}\right)\right).
\]
For example, when $p=2$ one has $\delta(\mathcal{N}_2)\ge\delta(\mathcal{T}_2)\approx 0.451699$.

It is not clear how Theorem~\ref{thm:sequence} can be efficiently used to improve this trivial lower bound
for $\delta(\mathcal{N}_p)$.
In the case of $p=2$ we sketch how to obtain a small improvement by considering certain nontrivial elements of $\mathcal{N}_2$
exhibited in~\cite{Mat:nonsing-der-2}.
It was shown there that $(2^{st}-1)/(2^s-1)$ belongs to $\mathcal{N}_2$, for $s\ge 1$ and $t\ge 3$.
This is a consequence of Theorem~\ref{thm:sequence} for $t\ge 4$, but not for $t=3$,
in which case the conclusion follows from a direct calculation given in~\cite[Proposition~2.2]{Mat:nonsing-der-2}.
Note that, by taking $s=1$, the numbers considered here include all numbers of the form $2^t-1$, except $3$.
Denoting by $\mathcal{S}$ the set of positive integers which are either multiples of $3$ or
of some number of the form $(2^{st}-1)/(2^s-1)$, with $s\ge 1$ and $t\ge 3$, we have
$\mathcal{T}_2\subseteq\mathcal{S}\subseteq\mathcal{N}_2$.
According to~\cite[Proposition~3.4]{Mat:nonsing-der-2} and the following comments,
$\mathcal{S}$ equals the set of multiples of numbers in the set
\[
\{3\}
\cup
\left\{\left.
\frac{2^{2^{a+2}}-1}{2^{2^a}-1}\,
\right\vert\,
a\ge 1
\right\}
\cup
\left\{\left.
\frac{2^{r^{b+1}}-1}{2^{r^b}-1}\,
\right\vert\,
r\text{ odd prime},\
b\ge 0
\right\}.
\]
Using the fact that the numbers of this set are mostly pairwise coprime,
with some exceptions detailed in~\cite[Proposition~3.4]{Mat:nonsing-der-2},
and performing some numerical calculations we have found that
$\delta(\mathcal{N}_2)\ge\delta(\mathcal{S})\approx 0.465673$.
By adding to the set $\mathcal{S}$ the multiples of the elements of $\mathcal{N}_2$
listed in Subsection~\ref{subsec:calculations} (which are the minimal elements of $\mathcal{N}_2$
not exceeding $200000$) we can still improve the lower bound for $\delta(\mathcal{N}_2)$ a little bit and obtain that
$\delta(\mathcal{N}_2)\ge 0.465926$.

\subsection{Relative size of elements of $\mathcal{N}_p$}

Our main result, especially in the slightly weaker but simpler form of Corollary~\ref{cor:n_large},
suggests that it may be interesting to introduce a relative measure of the size of a divisor $n$ of $p^k-1$.
There are several good candidates for this quantity, all close to $\log_p(n)/k$ for $n$ large, where
$k$ is taken as small as possible, hence $k=\ord_n(p)$, the (multiplicative) order of $p$ modulo $n$.
For the present discussion we select the following.
For an integer $n$ prime to $p$, we define its
{\em relative size} with respect to $p$ as the quantity $\log(n)/\log(p^{\ord_n(p)}-1)$.
With this definition, numbers of the form $p^k-1$ have relative size $1$.
Now we discuss what we know about the relative size of elements of $\mathcal{N}_p$.

The sufficient condition for $n\in\mathcal{N}_p$ given in Corollary~\ref{cor:n_large}
can be roughly read as the relative size of $n$ being slightly larger than $1-1/(2p)$.
However, the distinguished element $(p^p-1)/(p-1)$ of $\mathcal{N}_p$ has relative size close to $1-1/p$.
When $p=2$, Corollary~\ref{cor:n_large} says precisely that any odd number
of relative size at least $3/4$ belongs to $\mathcal{N}_2$.
The element $(2^{3s}-1)/(2^s-1)$ of $\mathcal{N}_2$ has relative size close to $2/3$, for $s$ large.
Several of the minimal elements of $\mathcal{N}_2$ listed in Subsection~\ref{subsec:calculations}
have even smaller relative size, the smallest, close to $0.433052$, being attained by $121369$.

Note that $\mathcal{N}_p$ contains elements of arbitrarily small relative size,
simply because it is closed with respect to taking multiples.
A simple way to see that is as follows.
For a prescribed element $n$ of $\mathcal{N}_p$ consider a prime $r$, different from $p$.
Then $nr$ belongs to $\mathcal{N}_p$ and has relative size
$\log(nr)/\log(p^{\ord_{nr}(p)}-1)\le
2\log(nr)/(\ord_r(p)\cdot\log p)$.
By an appropriate choice of $r$ the latter quantity can be made arbitrarily small, because
the function $\ord_r(p)/\log r$ of the prime $r$ is unbounded, as is easy to see
(in fact, much stronger statements hold, see~\cite{ErdMur}).
This argument, however, does not answer the question of whether {\em minimal} elements of $\mathcal{N}_p$
can have arbitrarily small relative size.

\section{The number of solutions of a certain system of equations}\label{sec:system}

The following proof of Theorem~\ref{thm:sequence} depends on Lemma~\ref{lemma:system}, which we postpone for clarity.

\begin{proof}[Proof of Theorem~\ref{thm:sequence}]
Let $N$ be the number of solutions over $\F_q$ of the system of equations
\begin{equation*}
\left\{
\begin{array}{rcl}
y_1^d &=& x
\\
y_2^d &=& x+1
\\
&\vdots& \\
y_r^d &=& x+r-1
\end{array}
\right.
\end{equation*}
An element $\xi$ of $\F_q$ such that
$\xi,\xi+1,\ldots,\xi+r-1$ are $d$th powers in $\F_q$
corresponds to $d^r$ distinct solutions of the system if none of the $\xi,\xi+1,\ldots,\xi+r-1$ equals zero,
and to $d^{r-1}$ solutions otherwise.
Since altogether these account for all solutions of the system, we have $N=d^rM+d^{r-1}M_0$,
and the desired inequality follows from Lemma~\ref{lemma:system}.
\end{proof}

\begin{lemma}\label{lemma:system}
Let $d$ be a divisor of $q-1$ and let $0<r\le p$.
Then the number $N$ of solutions over $\F_q$ of the system of equations
\begin{equation}\label{eq:system1}
\left\{
\begin{array}{rcl}
y_1^d &=& x
\\
y_2^d &=& x+1
\\
&\vdots& \\
y_r^d &=& x+r-1
\end{array}
\right.
\end{equation}
satisfies
\[
|N+d^{r-1}-q-1|\le\bigl((rd-r-d-1)d^{r-1}+2\bigr)\sqrt{q}.
\]
\end{lemma}

\begin{proof}
Let $\chi$ be a multiplicative character of $\F_q$ of (exact) order $d$.
Then all characters of order dividing $d$ are given by the powers $\chi^i$, for $i=0,\ldots,d-1$.
For each $j=1,\ldots,r$, and for any given $\xi\in\F_q$, the sum
\[
\sum_{i=0}^{d-1}\chi^i(\xi+j-1)
=
\sum_{i=0}^{d-1}\chi\bigl((\xi+j-1)^i\bigr)
\]
(reading $0^0$ as $1$ when it occurs as the argument of $\chi$)
equals the number of solutions of $y_j^d=\xi+j-1$.
Therefore, the product of all these quantities equals the number of solutions of the system
having $x=\xi$.
Consequently, the total number of solutions of system~\eqref{eq:system1} is given by
\begin{equation}\label{eq:N}
N=
\sum_{i_1=0}^{d-1}\cdots\sum_{i_r=0}^{d-1}
\sum_{\xi\in\F_q}\chi\bigl(\xi^{i_1}(\xi+1)^{i_2}\cdots(\xi+r-1)^{i_r}\bigr).
\end{equation}

It remains to evaluate or bound the character sum
$\sum_{\xi\in\F_q}\chi\bigl(\xi^{i_1}(\xi+1)^{i_2}\cdots(\xi+r-1)^{i_r}\bigr)$,
depending on the $r$tuple
$(i_1,\ldots,i_r)$.
The sum takes the value $q$ for the $r$tuple $(0,\ldots,0)$.
This case aside, the polynomial $z^{i_1}(z+1)^{i_2}\cdots(z+i-1)^{i_r}$
is never a $d$th power in $\F_p[z]$.
Therefore, Weil's bound for character sums~\cite[Theorem~5.41]{LN} applies and yields that
\begin{equation}\label{eq:Weil}
\biggl|
\sum_{\xi\in\F_q}\chi\bigl(\xi^{i_1}(\xi+1)^{i_2}\cdots(\xi+r-1)^{i_r}\bigr)
\biggr|
\le
\bigl(w(i_1,\ldots,i_r)-1\bigr)\sqrt{q},
\end{equation}
where $w(i_1,\ldots,i_r)$ is the number of distinct roots in $\F_q$ of the polynomial
$z^{i_1}(z+1)^{i_2}\cdots(z+i-1)^{i_r}$.
Clearly, $w(i_1,\ldots,i_r)$ equals the number of nonzero entries in the $r$tuple $(i_1,\ldots,i_r)$.
Adding together all character sums corresponding to the $r$tuples different from $(0,\ldots,0)$,
and using the triangle inequality, we obtain that $|N-q|$ does not exceed $\sqrt{q}$ times
the integer obtained by subtracting from $(d^{r}-1)(r-1)$
the total number of zero entries appearing in the collection of nonzero $r$tuples.
The total number of those zeroes equals $rd^{r-1}-r$, because
zero occurs as many times as any other integer $1,\ldots,d-1$
in the whole set of $r$tuples including $(0,\ldots,0)$.
We conclude that
\begin{equation}\label{eq:weak}
|N-q|\le\bigl((dr-r-d)d^{r-1}+1\bigr)\sqrt{q}.
\end{equation}

This inequality is close to our goal, but can still be improved a little (see Remark~\ref{rem:weak}).
The number of $r$tuples
$(i_1,\ldots,i_r)\neq(0,\ldots,0)$
such that $i_1+\cdots+i_r\equiv 0\pmod{d}$
is $d^{r-1}-1$.
Consider any one of them.
Then at least one of the entries $i_j$ is positive, say $i_1$ without loss of generality.
Since $\chi(\xi^d)=1$ for $\xi\in\F_q^\ast$ and $\chi(0)=0$, we have
\begin{align*}
\sum_{\xi\in\F_q}\chi\bigl(\xi^{i_1}(\xi+1)^{i_2}\cdots(\xi+r-1)^{i_r}\bigr)
&=
\sum_{\xi\in\F_q^\ast}\chi\bigl(\xi^{-i_2-i_3-\cdots-i_r}(\xi+1)^{i_2}\cdots(\xi+r-1)^{i_r}\bigr)
\\&=
\sum_{\xi\in\F_q^\ast}\chi\bigl((1+\xi^{-1})^{i_2}\cdots(1+(r-1)\xi^{-1})^{i_r}\bigr)
\\&=
\sum_{\eta\in\F_q^\ast}\chi\bigl((1+\eta)^{i_2}\cdots(1+(r-1)\eta)^{i_r}\bigr)
\\&=
-1+
\sum_{\eta\in\F_q}\chi\bigl((1+\eta)^{i_2}\cdots(1+(r-1)\eta)^{i_r}\bigr)
\end{align*}
The polynomial
$(1+z)^{i_2}\cdots(1+(r-1)z)^{i_r}$,
which provides the argument for $\chi$ in the last character sum,
has exactly $w(i_1,\ldots,i_r)-1$ distinct roots,
that is, one less than the polynomial corresponding to the original sum.
Therefore, for the character sums under present consideration inequality~\eqref{eq:Weil} can be strengthened to
\begin{equation*}
\biggl|
1+\sum_{\xi\in\F_q}\chi\bigl(\xi^{i_1}(\xi+1)^{i_2}\cdots(\xi+r-1)^{i_r}\bigr)
\biggr|
\le
\bigl(w(i_1,\ldots,i_r)-2\bigr)\sqrt{q}.
\end{equation*}
It follows that the coefficient of $\sqrt{q}$ in inequality~\eqref{eq:weak} can be decreased by $1$
for each of those $d^{r-1}-1$ character sums considered here, provided we increase $N-q$ by a constant term $1$ each time.
The desired inequality now follows.
\end{proof}

\begin{rem}\label{rem:weak}
The estimate for $M$ given in~\cite[Exercise~(5.66)]{LN}
(for a more general question, but that greater generality is inessential) is
\[
\left|M-\frac{q}{d^r}\right|
\le
\left(r-1-\frac{r}{d}+\frac{1}{d^r}\right)
\sqrt{q}+\frac{r}{d},
\]
and hence has the coefficient of $\sqrt{q}$ about $1/d$ larger than the estimate given in Theorem~\ref{thm:sequence}.
Since the values of $d$ of present interest to us may be much smaller than $\sqrt{q}$
(namely, roughly of the size of $q^{1/2p}$, see Corollary~\ref{cor:n_large}), this makes a significant difference
(in Theorem~\ref{thm:n_large}).
The larger coefficient of $\sqrt{q}$ given in~\cite[Exercise~(5.66)]{LN}
results from being content with inequality~\eqref{eq:weak}
in the proof of Lemma~\ref{lemma:system} (and thus, essentially, disregarding the effect of points at infinity).
For example, when $r=2$ it yields the weaker bound $|N-q|\le(d-1)^2\sqrt{q}$ rather than Weil's bound
$|N-d-q-1|\le(d-1)(d-2)\sqrt{q}$ for the Fermat curve $y_2^d-y_1^d=1$.
\end{rem}

\begin{rem}\label{rem:sophisticated}
The inequality proved in Lemma~\ref{lemma:system} is exactly Weil's bound
$|\bar N-q-1|\le 2g\sqrt{q}$ for the number $\bar N$ of $\F_q$-rational projective points
of the curve in the projective space $\mathbb{P}^{r+1}$ given by the system~\eqref{eq:system1} in affine coordinates.
In fact, the only singularity of the curve represented by~\eqref{eq:system1} occurs at its point at infinity,
which has multiplicity $d^{r-1}$.
An efficient way to compute the genus $g$ is to consider the nonsingular curve in $\mathbb{P}^r$,
which is birationally equivalent to~\eqref{eq:system1} via a projection,
given in affine coordinates by
\[
\left\{
\begin{array}{rcl}
y_2^d &=& y_1^d+1
\\
y_3^d &=& y_2^d+1
\\
&\vdots& \\
y_{r}^d &=& y_{r-1}^d+1
\end{array}
\right.
\]
Because this curve is (nonsingular and) a complete intersection of hypersurfaces,
one can compute its genus by means of the Adjunction Formula
and its iterates (see~\cite[V, Proposition~1.5]{Har}).
For a complete intersection of $s$ hypersurfaces of degrees $d_1,\ldots,d_s$, the Adjunction Formula reads
$2g-2=\bigl((\sum_i d_i)-s-2\bigr)\prod_i d_i$.
Since the curve under consideration is a complete intersection of $r-1$ hypersurfaces of degree $d$,
the Adjunction Formula gives
$2g-2=(rd-r-d-1)d^{r-1}$, as desired.
This argument gives a shorter but less elementary proof of Lemma~\ref{lemma:system}
based on Weil's bound $|\bar N-q-1|\le 2g\sqrt{q}$.
\end{rem}

\bibliography{References}

\end{document}